\documentclass[12pt,numbers]{elsarticle}
\journal{}

\makeatletter
\def\ps@pprintTitle{%
 \let\@oddhead\@empty
 \let\@evenhead\@empty
 \def\@oddfoot{\hfill\thepage}%
 \let\@evenfoot\@oddfoot}
\makeatother

\usepackage{amsmath,amssymb,amsfonts,amsthm,graphicx}
\usepackage[bookmarksnumbered,colorlinks=true]{hyperref}
\usepackage[labelfont=bf]{caption}

\providecommand{\doi}[1]{\href{https://doi.org/#1}{DOI:#1}}
\usepackage{xurl} 
\renewcommand{\doi}[1]{%
 \href{https://doi.org/#1}{\nolinkurl{DOI:#1}}%
}

\usepackage{dsfont} 
\usepackage{enumitem} 
\usepackage{mathtools} 
\usepackage{appendix} 
\usepackage{geometry} 
\usepackage{comment}
\numberwithin{equation}{section}
\numberwithin{table}{section}
\numberwithin{figure}{section}

\theoremstyle{plain}
\newtheorem{theorem}{Theorem}[section]

\newtheorem{lemma}[theorem]{Lemma}
\newtheorem{corollary}[theorem]{Corollary}

\theoremstyle{definition}
\newtheorem{remark}{Remark}[section]

\newcommand{\R}{\mathbb{R}}
\newcommand{\EE}{\mathsf{E}} 

\newcommand{\rd}{\mathrm{d}}

\newcommand{\leqdef}{\vcentcolon=}

\allowdisplaybreaks

\begin{document}

\begin{frontmatter}

\title{Log-convexity and log-concavity of noncentral gamma sums and differences}

\author[a2]{Robert E.~Gaunt}
\author[a3]{Fr\'ed\'eric Ouimet\corref{mycorrespondingauthor}}

\address[a2]{The University of Manchester, Manchester, M13 9PL, UK}
\address[a3]{Universit\'e du Qu\'ebec \`a Trois-Rivi\`eres, Trois-Rivi\`eres, QC G8Z 4M3, Canada\vspace{-5mm}}

\cortext[mycorrespondingauthor]{Corresponding author. Email address: frederic.ouimet2@uqtr.ca}

\begin{abstract}
We study log-convexity and log-concavity of densities obtained from sums and differences of two independent noncentral gamma random variables. We give a complete classification of one-sided log-convexity for noncentral gamma differences, a complete log-convexity classification for sums of two independent central gamma random variables, and sharp log-concavity criteria for central differences and for common-scale sums. As special cases, we deduce a log-convexity classification for the density of the product of two correlated normal random variables with arbitrary means and variances, and log-convexity and log-concavity classifications for the densities of the variance-gamma and McKay Type I distributions.
\end{abstract}


\end{frontmatter}

\section{Introduction}\label{sec:intro}

For $a,b > 0$ and $\lambda \geq 0$, we follow \cite{m93} and say that $G\sim\Gamma(a,b,\lambda)$ follows the noncentral gamma distribution if its characteristic function is given by
\begin{equation}\label{eq:noncentral.gamma.cf}
\EE[e^{itG}] = (1 - ibt)^{-a}\exp\left(\frac{ibt\lambda}{1 - ibt}\right), \qquad t\in\R.
\end{equation}
In the central case $\lambda = 0$, we recover the gamma distribution with probability density function
\begin{equation}\label{eq:central.gamma.density}
f_{a,b,0}(x) = \frac{x^{a - 1}e^{-x/b}}{b^a\Gamma(a)}, \qquad x > 0.
\end{equation}
For $\lambda > 0$, it was noted by \cite{GauntSutcliffe2026NCGamma} that the density is given by
\begin{equation}\label{eq:noncentral.gamma.density}
f_{a,b,\lambda}(x) = \frac{1}{b}e^{-\lambda - x/b}\left(\frac{x}{b\lambda}\right)^{(a - 1)/2}I_{a - 1}\left(2\sqrt{\frac{\lambda x}{b}}\right), \qquad x > 0,
\end{equation}
where $I_{\eta}$ is the modified Bessel function of the first kind; see, e.g., \cite[Chapter 10]{MR2723248}. Equivalently, by the Poisson mixture representation \cite{kb96,p49}, or directly from the defining power series for $I_{\eta}$ \citep[Eq.~10.25.2]{MR2723248},
\[
f_{a,b,\lambda}(x) = e^{-\lambda}\sum_{k=0}^{\infty}\frac{\lambda^k}{k!}\frac{x^{a + k - 1}e^{-x/b}}{b^{a + k}\Gamma(a + k)}, \qquad x > 0.
\]

The sum and difference of two independent noncentral gamma random variables, which are referred to by \cite{GauntSutcliffe2026NCGamma} as the noncentral gamma sum and difference distributions, arise in numerous settings. These distributions can equivalently be expressed as linear combinations of two independent noncentral chi-square random variables, which, in the integer-degrees-of-freedom case, can themselves be expressed as quadratic forms in normal random variables; see \cite{cl05,GauntSutcliffe2026NCGamma} for an overview of the extensive literature concerning such distributions. A review of the central gamma difference distribution is also given by \cite{k15}. Moreover, the sum and difference of two independent noncentral gamma random variables contain several important distributions as special cases, including the product of two correlated normal random variables \cite{craig,g26,wb32}, the variance-gamma distribution \cite{FischerGauntSarantsev2025VGReview} (also known as the generalized Laplace distribution \cite{kkp01}, the McKay Type II distribution \cite{HolmAlouini2004}, and the Bessel function distribution \cite{m32}), and the McKay Type I distribution \cite{HolmAlouini2004,m32} (also referred to as the McKay $I_\nu$ Bessel distribution \cite{pog1}).

Let $G_i\sim\Gamma(a_i,b_i,\lambda_i)$, $i = 1,2$, be independent noncentral gamma random variables, where $a_i,b_i > 0$ and $\lambda_i \geq 0$. The purpose of this paper is to determine sharp log-convexity and log-concavity criteria for the densities of the difference $D = G_1-G_2$ and the sum $S = G_1 + G_2$, including complete classifications in the central, one-sided and common-scale cases. We thus make a contribution to the body of literature on log-convexity and log-concavity of probability distributions, a subject of interest in many research domains such as economics \cite{a98}, reliability theory \cite{ss87} and queuing theory \cite{fw98}; see \cite{bb05} for basic theory and for a tabulation of standard probability distributions as having log-convex or log-concave densities.

Our first main result, Theorem \ref{thm:NCG.difference.shape}, is a complete one-sided log-convexity classification for the difference $D = G_1 - G_2$. We prove that the density of $D$ is log-convex on $(0,\infty)$ if and only if the first summand is central and has shape at most one, namely $\lambda_1 = 0$ and $a_1 \leq 1$. Similarly, it is log-convex on $(-\infty,0)$ if and only if $\lambda_2 = 0$ and $a_2 \leq 1$. The proof is based on a fixed-limit convolution formula for the density on each half-line. In the central case, the relevant kernel is a log-convex power times an exponential factor, whereas in the noncentral case the tail contains a positive square-root contribution in the logarithm, which rules out log-convexity by a simple convexity obstruction.

We also establish the corresponding results for sums in Theorem \ref{thm:NCG.sum.shape}. We prove that the density of $S = G_1 + G_2$ is log-convex on $(0,\infty)$ in the central case if and only if $a_1 + a_2 \leq 1$. We also prove log-concavity whenever both shapes are at least one, using the preservation of log-concavity under convolution; see \cite[Theorem~7]{Prekopa1973}. In the common-scale case $b_1 = b_2$, the sum is again noncentral gamma (see \cite[Remark 2.3]{GauntSutcliffe2026NCGamma}), and from shape results for the noncentral gamma distribution (which we record in Lemma \ref{prop:one.NCG.shape}) we obtain a complete log-concavity classification and a complete log-convexity classification for the sum.

Our results for the noncentral gamma difference and sum distributions unify several distributional shape questions. Indeed, as special cases, we obtain a one-sided log-convexity classification for the density of the product of two correlated normal random variables with arbitrary means and variances, and log-convexity and log-concavity classifications for the densities of the variance-gamma and McKay Type I distributions; see Corollaries \ref{cor:product.normal.classification}, \ref{cor:VG.shape} and \ref{cor:McKay.Type.I}, respectively. To the best of our knowledge, the product-normal and McKay Type I classifications in Corollaries~\ref{cor:product.normal.classification} and \ref{cor:McKay.Type.I} are new, while the variance-gamma classification in Corollary~\ref{cor:VG.shape} is included for completeness and also follows from known shape results for generalized hyperbolic densities \cite[Corollary~1]{YuNVM2011}. All results are stated in Section \ref{sec:results}, and all proofs are given in Section \ref{sec:proofs}.

\section{Results}\label{sec:results}

We begin by stating our conventions regarding log-convexity and log-concavity. A positive function $h:I\to(0,\infty)$ on an interval $I\subseteq\R$ is called log-convex if
\[
h(tx + (1 - t)y) \leq h(x)^t h(y)^{1 - t}, \qquad x,y\in I, \quad t\in[0,1].
\]
It is called log-concave if the reverse inequality holds. Equivalently, $\log(h)$ is convex in the log-convex case and concave in the log-concave case. In the present setting, these notions will always be understood on the interior of the support under consideration.

The following elementary convexity and concavity facts will be used throughout.

\begin{lemma}[Convexity and concavity facts]\label{lem:log.convexity.facts}
Let $(E,\mathcal{E},\eta)$ be a positive measure space. Suppose that, for every $x\in I$, the map $y\mapsto h_y(x)$ is measurable. If $h_y:I\to(0,\infty)$ is log-convex for every $y\in E$ and
\[
H(x) \leqdef \int_E h_y(x)\,\eta(\rd y)
\]
is finite and strictly positive for every $x\in I$, then $H$ is log-convex on $I$. Also, if $G$ is convex on $[x_0,\infty)$ and $G(x)/x\to0$ as $x\to\infty$, then $G$ cannot satisfy $\limsup_{x\to\infty}G(x) = \infty$. Finally, the convolution of two log-concave densities, extended by zero outside their supports, is log-concave.
\end{lemma}

The starting point of our analysis is the shape of a single noncentral gamma density. The following lemma, which follows from basic shape results for the gamma distribution, a known shape classification for the noncentral chi-square density \cite[Theorem 1, i, ii]{y11}, and a simple tail argument in the noncentral log-convexity case, is the local building block for the convolution arguments given in the proofs of Theorems \ref{thm:NCG.difference.shape} and \ref{thm:NCG.sum.shape}, and the lemma also explains why the log-convex case necessarily requires centrality.

\begin{lemma}[Shape of one noncentral gamma density]\label{prop:one.NCG.shape}
Let $a,b > 0$ and $\lambda \geq 0$. The density $f_{a,b,\lambda}$ of $G\sim\Gamma(a,b,\lambda)$ is log-convex on $(0,\infty)$ if and only if $\lambda = 0$ and $a \leq 1$. The density $f_{a,b,\lambda}$ is log-concave on $(0,\infty)$ if and only if $a \geq 1$.
\end{lemma}

We now pass from one density to the difference of two independent noncentral gamma random variables. The density on each open half-line has a fixed-limit convolution representation, so one-sided log-convexity is controlled by the gamma density that appears with the shifted argument. The following theorem shows that this control is exact, and that in the central case the natural log-concavity condition is also necessary.

\begin{theorem}[Shape of noncentral gamma differences]\label{thm:NCG.difference.shape}
Let $G_i\sim\Gamma(a_i,b_i,\lambda_i)$, $i = 1,2$, be independent, where $a_i,b_i > 0$ and $\lambda_i \geq 0$. Let $p_{-}$ be the density of
\[
D \leqdef G_1 - G_2.
\]
Then the following assertions hold.
\begin{enumerate}[label = {\rm(\roman*)}]\setlength\itemsep{0em}
\item The density $p_{-}$ is log-convex on $(0,\infty)$ if and only if $\lambda_1 = 0$ and $a_1 \leq 1$.
\item The density $p_{-}$ is log-convex on $(-\infty,0)$ if and only if $\lambda_2 = 0$ and $a_2 \leq 1$.
\item If $a_1 \geq 1$ and $a_2 \geq 1$, then $p_{-}$ is log-concave on $\R$.
\item If $\lambda_1 = \lambda_2 = 0$, then $p_{-}$ is log-concave on $\R$ if and only if $a_1 \geq 1$ and $a_2 \geq 1$.
\end{enumerate}
\end{theorem}

\newpage
Our corresponding result for the sum of two independent noncentral gamma random variables is given in the following theorem.

\begin{theorem}[Shape of noncentral gamma sums]\label{thm:NCG.sum.shape}
Let $G_i\sim\Gamma(a_i,b_i,\lambda_i)$, $i = 1,2$, be independent, where $a_i,b_i > 0$ and $\lambda_i \geq 0$. Let $p_{+}$ be the density of
\[
S \leqdef G_1 + G_2.
\]
Then the following assertions hold.
\begin{enumerate}[label = {\rm(\roman*)}]\setlength\itemsep{0em}
\item If $\lambda_1 = \lambda_2 = 0$, then $p_{+}$ is log-convex on $(0,\infty)$ if and only if $a_1 + a_2 \leq 1$.
\item If $a_1 \geq 1$ and $a_2 \geq 1$, then $p_{+}$ is log-concave on $(0,\infty)$.
\item If $b_1 = b_2$, then $p_{+}$ is log-concave on $(0,\infty)$ if and only if $a_1 + a_2 \geq 1$, and $p_{+}$ is log-convex on $(0,\infty)$ if and only if $\lambda_1 = \lambda_2 = 0$ and $a_1 + a_2 \leq 1$.
\item Suppose $\lambda_1 = \lambda_2 = 0$ and write $a \leqdef a_1 = a_2$. If $b_1 \neq b_2$, then $p_{+}$ is log-concave on $(0,\infty)$ if and only if $a \geq 1$.
\end{enumerate}
\end{theorem}

In the following corollary, we translate the noncentral gamma statements into the more familiar language of weighted noncentral chi-square variables.

\begin{corollary}[Linear combinations of two noncentral chi-square variables]\label{cor:chisquare.linear.combinations}
Let $V_i\sim \chi^2_{\nu_i}(\delta_i)$, $i = 1,2$, be independent, where $\nu_i > 0$ and $\delta_i \geq 0$, and let $w_i > 0$. Let $r_{-}$ and $r_{+}$ be the densities of $w_1V_1 - w_2V_2$ and $w_1V_1 + w_2V_2$, respectively. Then the following assertions hold.
\begin{enumerate}[label = {\rm(\alph*)}]\setlength\itemsep{0em}
\item The density $r_{-}$ is log-convex on $(0,\infty)$ if and only if $\delta_1 = 0$ and $\nu_1 \leq 2$.
\item The density $r_{-}$ is log-convex on $(-\infty,0)$ if and only if $\delta_2 = 0$ and $\nu_2 \leq 2$.
\item If $\delta_1 = \delta_2 = 0$, then $r_{+}$ is log-convex on $(0,\infty)$ if and only if $\nu_1 + \nu_2 \leq 2$.
\item If $\nu_1 \geq 2$ and $\nu_2 \geq 2$, then $r_{-}$ is log-concave on $\R$, and $r_{+}$ is log-concave on $(0,\infty)$.
\item If $\delta_1 = \delta_2 = 0$, then $r_{-}$ is log-concave on $\R$ if and only if $\nu_1 \geq 2$ and $\nu_2 \geq 2$.
\item If $w_1 = w_2$, then $r_{+}$ is log-concave on $(0,\infty)$ if and only if $\nu_1 + \nu_2 \geq 2$, and $r_{+}$ is log-convex on $(0,\infty)$ if and only if $\delta_1 = \delta_2 = 0$ and $\nu_1 + \nu_2 \leq 2$.
\end{enumerate}
\end{corollary}

Let $(X,Y)$ be a bivariate normal random vector with mean vector $(\mu_X,\mu_Y)\in \R^2$, positive variances $(\sigma_X^2,\sigma_Y^2)$, and correlation coefficient $\rho\in(-1,1)$. The following corollary provides a complete one-sided log-convexity classification for the density of the product $Z = XY$. The corollary is a consequence of Corollary \ref{cor:chisquare.linear.combinations} and the following connection between the distribution of the product $Z$ and a weighted difference of two independent noncentral chi-square variates. By Theorem 2.1 of \cite{g26} (see also \cite{j26} for an alternative proof), we have
\begin{align}\label{conec}
Z=_d\frac{\sigma_X\sigma_Y}{2}(1 + \rho)V_1-\frac{\sigma_X\sigma_Y}{2}(1 - \rho)V_2,
\end{align}
where $V_1\sim \chi_1^2(\lambda_{+})$ and $V_2\sim \chi_1^2(\lambda_{-})$ are independent noncentral chi-square random variables with one degree of freedom and noncentrality parameters
\begin{align*}
\lambda_{+} = \frac{(\mu_X/\sigma_X + \mu_Y/\sigma_Y)^2}{2(1 + \rho)}, \qquad \lambda_{-} = \frac{(\mu_X/\sigma_X - \mu_Y/\sigma_Y)^2}{2(1 - \rho)}.
\end{align*}

\begin{corollary}[Log-convexity classification for the product of two correlated normal random variables]\label{cor:product.normal.classification}
The density $p_Z$ of the product $Z = XY$ is log-convex on $(0,\infty)$ if and only if
\[
\frac{\mu_X}{\sigma_X} + \frac{\mu_Y}{\sigma_Y} = 0.
\]
The density $p_Z$ is log-convex on $(-\infty,0)$ if and only if
\[
\frac{\mu_X}{\sigma_X} - \frac{\mu_Y}{\sigma_Y} = 0.
\]
Consequently, $p_Z$ is log-convex on both open half-lines if and only if $\mu_X = \mu_Y = 0$.
\end{corollary}

\begin{remark}
It is readily seen that the density $p_Z$ is not log-concave on either of the intervals $(-\infty,0)$ or $(0,\infty)$, since, for all parameter values, the density has a logarithmic singularity at the origin; see \cite[Proposition 2.1]{gz23}.
\end{remark}

For parameters $\nu > -1/2$, $\alpha > 0$, $\beta\in\R$, $|\beta| < \alpha$, and $\mu\in\R$, the variance-gamma distribution, denoted by $\mathrm{VG}(\nu,\alpha,\beta,\mu)$, has density
\begin{equation}\label{eq:variance.gamma.par}
p_{\mathrm{VG}}(x) = \frac{(\alpha^2 - \beta^2)^{\nu + 1/2}}{\sqrt{\pi}(2\alpha)^{\nu}\Gamma(\nu + 1/2)} e^{\beta(x - \mu)}|x - \mu|^{\nu} K_{\nu}(\alpha|x - \mu|), \qquad x\in\R\setminus\{\mu\};
\end{equation}
for this and other parameterizations of the variance-gamma distribution, see the review \cite{FischerGauntSarantsev2025VGReview}. Here $K_{\eta}$ is the modified Bessel function of the second kind (see \cite[Chapter 10]{MR2723248}). A variance-gamma random variable with location parameter $\mu = 0$ can be expressed as a difference of two independent central gamma variates (see \cite[Eq.~(20)]{FischerGauntSarantsev2025VGReview}), and since log-convexity and log-concavity are preserved under translations, the following corollary can be established by applying Theorem \ref{thm:NCG.difference.shape} and the one-sided tail argument used in its proof.

\begin{corollary}[Variance-gamma densities]\label{cor:VG.shape}
The variance-gamma density $p_{\mathrm{VG}}$ is log-convex on each of the intervals $(-\infty,\mu)$ and $(\mu,\infty)$ if and only if $-1/2 < \nu \leq 1/2$. The density is log-concave on each of these intervals if and only if $\nu \geq 1/2$. At $\nu = 1/2$, the density is log-affine on each of these intervals.
\end{corollary}

\begin{remark}
1. Corollary \ref{cor:VG.shape} asserts that the variance-gamma density is log-affine on the intervals $(-\infty,\mu)$ and $(\mu,\infty)$ when $\nu = 1/2$. As noted by \cite{kkp01}, this case corresponds to the asymmetric Laplace distribution with density
\[
p_{\mathrm{AL}}(x) = \frac{\alpha^2-\beta^2}{2\alpha}e^{\beta(x-\mu)-\alpha|x-\mu|}, \qquad x\in\R,
\]
which is immediately seen to be log-affine on $(-\infty,\mu)$ and $(\mu,\infty)$.

\vspace{3mm}

\noindent 2. When $\mu_X = \mu_Y = 0$, the distribution of the product of two correlated normal random variables $Z = XY$ is variance-gamma distributed with density
\[
p_Z(x) = \frac{1}{\pi\sigma_X\sigma_Y\sqrt{1-\rho^2}}\exp\bigg(\frac{\rho x}{\sigma_X\sigma_Y(1-\rho^2)}\bigg)K_0\bigg(\frac{|x|}{\sigma_X\sigma_Y(1-\rho^2)}\bigg), \quad x\in\R\setminus\{0\},
\]
(see \cite{gaunt prod,np16}). Thus, log-convexity of the product $Z$ on both open half-lines in the zero-mean case is consistent with the $\nu = 0$ case of Corollary~\ref{cor:VG.shape}. We note that if either $\mu_X$ or $\mu_Y$ is nonzero, the density $p_Z$ takes a rather complicated form, with the exact formulas of \cite{cui,MR4924268} expressed in terms of infinite series or integrals involving special functions. Therefore, the probabilistic representation of the product $Z$ as a weighted difference of two independent noncentral chi-square variates is crucial for establishing the log-convexity classification.
\end{remark}

Finally, we consider the McKay Type I distribution \cite{m32}, which, for $b > 0$, $c > 1$, $m > -1/2$, has density
\begin{equation}\label{mden}
p_{\mathrm{McKay}\, \mathrm{I}}(x) = \frac{\sqrt{\pi}(c^2-1)^{m + 1/2}}{2^m b^{m + 1}\Gamma(m + 1/2)}x^m e^{-cx/b}I_m\bigg(\frac{x}{b}\bigg), \qquad x > 0.
\end{equation}
It was shown by \cite[Theorem 3]{HolmAlouini2004} that a McKay Type I random variable can be expressed as the sum of two central gamma random variables. Thus, the following corollary can be deduced as a consequence of Theorem \ref{thm:NCG.sum.shape}.

\begin{corollary}[McKay Type I densities]\label{cor:McKay.Type.I}
The McKay Type I density $p_{\mathrm{McKay}\, \mathrm{I}}$ is log-convex on $(0,\infty)$ if and only if $m\leq0$. The density $p_{\mathrm{McKay}\, \mathrm{I}}$ is log-concave on $(0,\infty)$ if and only if $m \geq 1/2$.
\end{corollary}

\begin{remark}
For $0 < m < 1/2$, the McKay Type I density $p_{\mathrm{McKay}\, \mathrm{I}}$ is neither log-convex nor log-concave on $(0,\infty)$. In fact, for $0 < m < 1/2$, the density is log-concave near the origin and log-convex for sufficiently large $x$. Indeed, using the power series representation of the modified Bessel function of the first kind (see \citep[Eq.~10.25.2]{MR2723248}),
\begin{equation}\label{ismall}
I_{\eta}(x) = \sum_{k=0}^{\infty} \frac{1}{k!\Gamma(k + \eta + 1)}\bigg(\frac{x}{2}\bigg)^{2k + \eta} = \frac{x^{\eta}}{2^{\eta}\Gamma(\eta + 1)}\big(1 + \mathcal{O}(x^2)\big), \qquad x\downarrow0, \quad \eta > -1,
\end{equation}
and a straightforward asymptotic analysis, we obtain that, for $m > 0$,
\[
\big(\log\,p_{\mathrm{McKay}\, \mathrm{I}}\big)''(x)\sim -\frac{2m}{x^2}, \qquad x\downarrow0,
\]
so that the density $p_{\mathrm{McKay}\, \mathrm{I}}$ is log-concave near the origin for $0 < m < 1/2$. Also, applying the asymptotic expansion
\begin{equation}\label{ilarge}
I_{\eta}(x) = \frac{e^x}{\sqrt{2\pi x}}\big(1 + \mathcal{O}(x^{-1})\big), \qquad x\rightarrow\infty,
\end{equation}
(see \cite[Eq.~10.30.4]{MR2723248}), we obtain that
\[
\big(\log\,p_{\mathrm{McKay}\, \mathrm{I}}\big)''(x)\sim -\frac{m-1/2}{x^2}, \qquad x\rightarrow\infty,
\]
and therefore the density $p_{\mathrm{McKay}\, \mathrm{I}}$ is log-convex for sufficiently large $x$ when $0 < m < 1/2$.
\end{remark}

\begin{remark}
Corollaries~\ref{cor:VG.shape} and \ref{cor:McKay.Type.I} may also be read as shape classifications for two classical products of modified Bessel functions. Since multiplication by an exponential factor and rescaling of the argument preserve log-convexity and log-concavity, Corollary~\ref{cor:VG.shape} implies that, for $\nu > -1/2$, the function $x\mapsto x^{\nu} K_{\nu}(x)$ is log-convex on $(0,\infty)$ if and only if $-1/2 < \nu \leq 1/2$, and is log-concave on $(0,\infty)$ if and only if $\nu \geq 1/2$. This classification is also contained in the log-convexity and log-concavity criteria for generalized hyperbolic densities given by \cite[Corollary~1]{YuNVM2011}. Similarly, Corollary~\ref{cor:McKay.Type.I} implies that, for $\nu > -1/2$, the function $x\mapsto x^{\nu} I_{\nu}(x)$ is log-convex on $(0,\infty)$ if and only if $-1/2 < \nu \leq 0$, and is log-concave on $(0,\infty)$ if and only if $\nu \geq 1/2$. The latter log-concavity assertion is consistent with known log-concavity results for modified Bessel functions of the first kind, such as the log-concavity of $x\mapsto \sqrt{x}I_{\nu}(x)$ for $\nu \geq 1/2$ proved in \cite[Theorem~1(c)]{BariczPonnusamyVuorinen2011}, which implies log-concavity of $x\mapsto x^\nu I_{\nu}(x)$ for $\nu \geq 1/2$.
\end{remark}

\section{Proofs}\label{sec:proofs}

\subsection{Proof of Lemma~\ref{lem:log.convexity.facts}}

For $x_1,x_2\in I$ and $t\in[0,1]$, log-convexity and H\"older's inequality give
\[
H(tx_1 + (1 - t)x_2) \leq \int_E h_y(x_1)^t h_y(x_2)^{1 - t}\,\eta(\rd y) \leq H(x_1)^tH(x_2)^{1 - t}.
\]
This proves the first claim. For the second claim, suppose that $G(x_2) > G(x_1)$ for some $x_0 \leq x_1 < x_2$. Convexity gives, for every $x \geq x_2$,
\[
G(x) \geq G(x_2) + \frac{G(x_2) - G(x_1)}{x_2 - x_1}(x - x_2),
\]
which implies $\liminf_{x\to\infty}G(x)/x > 0$. This contradicts $G(x)/x\to0$. Hence no such pair $x_1,x_2$ can occur if $G(x)/x\to0$, and $G$ is bounded above by $G(x_0)$ on $[x_0,\infty)$. The convolution-preservation statement for log-concavity follows from \cite[Theorem~7]{Prekopa1973}.

\subsection{Proof of Lemma~\ref{prop:one.NCG.shape}}

When $\lambda = 0$, the assertion is the standard log-convexity and log-concavity classification for the gamma distribution; see, for example, \cite[Tables 1 and 2]{bb05}. Since a noncentral gamma random variable can be expressed as a scaled noncentral chi-square variate, for $\lambda > 0$, the log-concavity classification is an immediate consequence of the corresponding classification for the noncentral chi-square distribution given by \cite[Theorem 1, i]{y11}. It remains only to rule out log-convexity when $\lambda > 0$. The result \cite[Theorem 1, ii]{y11} rules this out when $0 < a < 1$, since the associated chi-square degrees of freedom are $\nu = 2a$, but it does not cover $a \geq 1$. The following tail argument applies for every $a > 0$. From \eqref{eq:noncentral.gamma.density} and \eqref{ilarge},
\[
f_{a,b,\lambda}(x) = Cx^{(2a - 3)/4}\exp\left(-\frac{x}{b} + 2\sqrt{\frac{\lambda x}{b}}\right)(1 + o(1)), \qquad x\to\infty,
\]
where $C \in(0,\infty)$. Hence $G(x)\leqdef \log\{f_{a,b,\lambda}(x)\} + x/b$ is unbounded above and satisfies $G(x)/x\to0$. If $f_{a,b,\lambda}$ were log-convex, then $G$ would be convex, contradicting Lemma~\ref{lem:log.convexity.facts}.

\subsection{Proof of Theorem~\ref{thm:NCG.difference.shape}}

For $x > 0$, the density of $D = G_1 - G_2$ is given by
\begin{equation}\label{eq:NCG.difference.positive.convolution}
p_{-}(x) = \int_0^{\infty} f_{a_1,b_1,\lambda_1}(x + y)f_{a_2,b_2,\lambda_2}(y)\,\rd y.
\end{equation}
Suppose first that $\lambda_1 = 0$ and $a_1 \leq 1$. By \eqref{eq:central.gamma.density},
\[
p_{-}(x) = \frac{e^{-x/b_1}}{b_1^{a_1}\Gamma(a_1)}\int_0^{\infty} (x + y)^{a_1 - 1}e^{-y/b_1}f_{a_2,b_2,\lambda_2}(y)\,\rd y.
\]
For every fixed $y \geq 0$, the function $x\mapsto (x + y)^{a_1 - 1}$ is log-convex on $(0,\infty)$, and the exponential factor has affine logarithm. Lemma~\ref{lem:log.convexity.facts} therefore gives log-convexity of $p_{-}$ on $(0,\infty)$.

Conversely, suppose that $p_{-}$ is log-convex on $(0,\infty)$. If $\lambda_1 > 0$, then recalling the expression \eqref{eq:noncentral.gamma.density} for the density of the noncentral gamma distribution and applying the asymptotic expansion \eqref{ilarge} inside \eqref{eq:NCG.difference.positive.convolution} yields
\begin{equation}\label{eq:NCG.difference.noncentral.tail}
p_{-}(x) = C_1x^{(2a_1 - 3)/4}\exp\left(-\frac{x}{b_1} + 2\sqrt{\frac{\lambda_1x}{b_1}}\right)(1 + o(1)), \qquad x\to\infty,
\end{equation}
where $C_1 \in(0,\infty)$. Indeed, from \eqref{ilarge} we have
\begin{align*}
\frac{f_{a_1,b_1,\lambda_1}(x + y)}{x^{(2a_1 - 3)/4}\exp(-x/b_1 + 2\sqrt{\lambda_1x/b_1})} \rightarrow C e^{-y/b_1}, \qquad x\rightarrow\infty,
\end{align*}
for some $C\in(0,\infty)$.
We note that the application of the asymptotic expansion \eqref{ilarge} inside \eqref{eq:NCG.difference.positive.convolution} is justified since the integrand has a uniform asymptotic expansion in the integration variable \cite{lopez}, which in our case is readily seen from the limiting forms \eqref{ismall} and \eqref{ilarge}.
It follows from \eqref{eq:NCG.difference.noncentral.tail} that $G(x) \leqdef \log\{p_{-}(x)\} + x/b_1$ is unbounded above and satisfies $G(x)/x\to0$. Since $G$ would be convex if $p_{-}$ were log-convex, this contradicts Lemma~\ref{lem:log.convexity.facts}. Therefore $\lambda_1 = 0$.

Under $\lambda_1 = 0$, for every $a_1 > 0$,
\begin{equation}\label{eq:NCG.difference.central.tail}
p_{-}(x) = C_2x^{a_1 - 1}e^{-x/b_1}(1 + o(1)), \qquad x\to\infty,
\end{equation}
where $C_2 \in(0,\infty)$. Indeed, after division by $x^{a_1 - 1}e^{-x/b_1}$, the integrand in \eqref{eq:NCG.difference.positive.convolution} converges pointwise to a positive constant times $e^{-y/b_1}f_{a_2,b_2,\lambda_2}(y)$, and it is bounded by a constant multiple of
\[
(1 + y)^{\max\{a_1 - 1,0\}}e^{-y/b_1}f_{a_2,b_2,\lambda_2}(y).
\]
Dominated convergence proves \eqref{eq:NCG.difference.central.tail}.
If $a_1 > 1$, then $G(x) \leqdef \log\{p_{-}(x)\} + x/b_1$ is again convex under the assumption of log-convexity, while \eqref{eq:NCG.difference.central.tail} shows that $G(x)\to\infty$ and $G(x)/x\to0$. This contradicts Lemma~\ref{lem:log.convexity.facts}. Hence $a_1 \leq 1$, proving part {\rm(i)}. Part {\rm(ii)} follows by applying part {\rm(i)} to $-D = G_2 - G_1$.

For part {\rm(iii)}, Lemma~\ref{prop:one.NCG.shape} shows that $f_{a_i,b_i,\lambda_i}$ is log-concave on $(0,\infty)$ when $a_i \geq 1$. Extending these densities by zero outside $(0,\infty)$ gives log-concave functions on $\R$. Since $p_{-}$ is the convolution of $f_{a_1,b_1,\lambda_1}$ and the reflection of $f_{a_2,b_2,\lambda_2}$, the convolution-preservation statement for log-concavity stated in Lemma~\ref{lem:log.convexity.facts} gives log-concavity of $p_{-}$ on $\R$.

It remains to prove the necessity in part {\rm(iv)}. Suppose that $\lambda_1 = \lambda_2 = 0$. If $a_1 < 1$, then \eqref{eq:NCG.difference.central.tail} gives
\[
\frac{p_{-}(2x)^2}{p_{-}(x)p_{-}(3x)}\to\left(\frac{4}{3}\right)^{a_1 - 1} < 1, \qquad x\to\infty.
\]
This contradicts the midpoint inequality required by log-concavity, for all sufficiently large $x$. If $a_2 < 1$, the same argument applied to $-D$ rules out log-concavity on $\R$. Therefore $a_1 \geq 1$ and $a_2 \geq 1$ are necessary. Their sufficiency is part {\rm(iii)}.

\begin{remark}
The asymptotic approximation \eqref{eq:NCG.difference.noncentral.tail} can also be recovered immediately from the asymptotic expansion for the density of the noncentral gamma distribution given in \cite[Theorem 3.4, part 1]{GauntSutcliffe2026NCGamma}. We have chosen to provide an efficient derivation of the limiting form \eqref{eq:NCG.difference.noncentral.tail} via the convolution formula \eqref{eq:NCG.difference.positive.convolution} to ensure that the proof of Theorem \ref{thm:NCG.difference.shape} is self-contained, up to basic convexity and concavity facts.
\end{remark}

\subsection{Proof of Theorem~\ref{thm:NCG.sum.shape}}

Assume first that $\lambda_1 = \lambda_2 = 0$. For $x > 0$,
\begin{equation}\label{intproof}
p_{+}(x) = C \int_0^x (x - y)^{a_1 - 1}y^{a_2 - 1}\exp\left(-\frac{x - y}{b_1} - \frac{y}{b_2}\right)\,\rd y,
\end{equation}
where $C \in(0,\infty)$. After the change of variable $y = xu$, we obtain
\[
p_{+}(x) = C e^{-x/b_1}\int_0^1 x^{a_1 + a_2 - 1}(1 - u)^{a_1 - 1}u^{a_2 - 1}\exp\left(-xu\left(\frac{1}{b_2} - \frac{1}{b_1}\right)\right)\,\rd u.
\]
If $a_1 + a_2 \leq 1$, then, for every fixed $u\in(0,1)$, the function of $x$ inside the integral is log-convex on $(0,\infty)$. Hence Lemma~\ref{lem:log.convexity.facts} gives log-convexity of $p_{+}$. Conversely,
\[
p_{+}(x)\sim C \frac{\Gamma(a_1) \Gamma(a_2)}{\Gamma(a_1 + a_2)} x^{a_1 + a_2 - 1}, \qquad x\downarrow0.
\]
If $a_1 + a_2 > 1$, then
\[
\frac{p_{+}(2x)^2}{p_{+}(x)p_{+}(3x)}\to\left(\frac{4}{3}\right)^{a_1 + a_2 - 1} > 1, \qquad x\downarrow0.
\]
This contradicts the midpoint inequality required by log-convexity, for all sufficiently small $x$. This proves part {\rm(i)}.

Part {\rm(ii)} follows from Lemma~\ref{prop:one.NCG.shape} and the preservation of log-concavity under convolution. Finally, if $b_1 = b_2 = b$, then the characteristic function formula \eqref{eq:noncentral.gamma.cf} gives
\[
G_1 + G_2\sim\Gamma(a_1 + a_2,b,\lambda_1 + \lambda_2).
\]
Part {\rm(iii)} now follows from Lemma~\ref{prop:one.NCG.shape}, because $\lambda_1 + \lambda_2 = 0$ if and only if $\lambda_1 = \lambda_2 = 0$.

Finally, we establish part {\rm(iv)}. Suppose now that $b_1 \neq b_2$ and $\lambda_1 = \lambda_2 = 0$, and let $a \leqdef a_1 = a_2$. If $a \geq 1$, then log-concavity follows from part {\rm(ii)}. Conversely, assume without loss of generality that $b_1 > b_2$. Put $\theta \leqdef 1/b_2 - 1/b_1 > 0$. From \eqref{intproof},
\[
p_{+}(x) = Ce^{-x/b_1}\int_0^x (x - y)^{a - 1}y^{a - 1}e^{-\theta y}\,\rd y.
\]
By splitting the integral over $(0,x/2)$ and $(x/2,x)$ and applying dominated convergence on the first interval, we obtain
\[
p_{+}(x) = C_1x^{a - 1}e^{-x/b_1}(1 + o(1)), \qquad x\to\infty,
\]
where $C_1 \in(0,\infty)$. If $a < 1$, then
\[
\frac{p_{+}(2x)^2}{p_{+}(x)p_{+}(3x)}\to\left(\frac{4}{3}\right)^{a - 1} < 1, \qquad x\to\infty.
\]
This contradicts the midpoint inequality required by log-concavity, for all sufficiently large $x$. This proves the claim.

\subsection{Proof of Corollary~\ref{cor:chisquare.linear.combinations}}

The distributional relation $w_iV_i\sim\Gamma(\nu_i/2,2w_i,\delta_i/2)$, $i = 1,2$, follows from \eqref{eq:noncentral.gamma.cf} and the standard characteristic function formula for the noncentral chi-square distribution; see also \cite[Remark~2.2 and Eq.~(2.3)]{GauntSutcliffe2026NCGamma}. Parts {\rm(a)} and {\rm(b)} are therefore exactly Theorem~\ref{thm:NCG.difference.shape}{\rm(i)} and {\rm(ii)}. In the central case $\delta_1 = \delta_2 = 0$, part {\rm(c)} is Theorem~\ref{thm:NCG.sum.shape}{\rm(i)}. Part {\rm(d)} follows from Theorem~\ref{thm:NCG.difference.shape}{\rm(iii)} and Theorem~\ref{thm:NCG.sum.shape}{\rm(ii)}, and part {\rm(e)} follows from Theorem~\ref{thm:NCG.difference.shape}{\rm(iv)}. Finally, part {\rm(f)} follows from Theorem~\ref{thm:NCG.sum.shape}{\rm(iii)}.

\subsection{Proof of Corollary~\ref{cor:product.normal.classification}}

We recall that the product $Z = XY$ can be expressed as a weighted difference of two independent noncentral chi-square variates through the distributional relation \eqref{conec}. Corollary \ref{cor:chisquare.linear.combinations}{\rm(i)} gives log-convexity on $(0,\infty)$ if and only if $\lambda_{+} = 0$, which is equivalent to $\mu_X/\sigma_X + \mu_Y/\sigma_Y = 0$. Similarly, Corollary~\ref{cor:chisquare.linear.combinations}{\rm(ii)} gives log-convexity on $(-\infty,0)$ if and only if $\lambda_{-} = 0$, which is equivalent to $\mu_X/\sigma_X - \mu_Y/\sigma_Y = 0$. Both conditions hold if and only if $\mu_X/\sigma_X + \mu_Y/\sigma_Y = \mu_X/\sigma_X - \mu_Y/\sigma_Y = 0$, which is equivalent to $\mu_X = \mu_Y = 0$.

\subsection{Proof of Corollary~\ref{cor:VG.shape}}

It suffices to take $\mu = 0$, because translations preserve log-convexity and log-concavity on the corresponding open half-lines. Put $\kappa \leqdef \nu + 1/2$, $b_1 \leqdef 1/(\alpha - \beta)$ and $b_2 \leqdef 1/(\alpha + \beta)$. Then \cite[Eq.~(20)]{FischerGauntSarantsev2025VGReview}, translated into the parameterization in \eqref{eq:variance.gamma.par}, gives the distributional relation
\[
G_1 - G_2\sim \mathrm{VG}(\nu,\alpha,\beta,0),
\]
where $G_1\sim\Gamma(\kappa,b_1,0)$ and $G_2\sim\Gamma(\kappa,b_2,0)$ are independent. Theorem~\ref{thm:NCG.difference.shape}{\rm(i)} and {\rm(ii)} give log-convexity on both open half-lines if and only if $\kappa \leq 1$, which is equivalent to $\nu \leq 1/2$. If $\kappa \geq 1$, then Theorem~\ref{thm:NCG.difference.shape}{\rm(iii)} gives log-concavity on $\R$, and hence on both open half-lines. Conversely, if $\kappa < 1$, then the one-sided tail argument in the proof of Theorem~\ref{thm:NCG.difference.shape}{\rm(iv)} rules out log-concavity on each open half-line. Thus log-concavity on each open half-line holds if and only if $\kappa \geq 1$, which is equivalent to $\nu \geq 1/2$. When $\nu = 1/2$, the log-convexity and log-concavity classifications meet, and the density is therefore log-affine on each open half-line.

\subsection{Proof of Corollary~\ref{cor:McKay.Type.I}}

Let $Y_1$ and $Y_2$ be independent $\Gamma(m + 1/2,1)$ random variables. Then, as noted by \cite[Eq.\ (3.44)]{gm21}, an application of \cite[Theorem 3]{HolmAlouini2004} implies that a McKay Type I random variable $Z_{m,c,b}$ with density \eqref{mden} satisfies the distributional relation
\[
Z_{m,c,b}=_d \frac{b}{c + 1}Y_1 + \frac{b}{c - 1}Y_2.
\]
The log-convexity assertion now follows from Theorem~\ref{thm:NCG.sum.shape}{\rm(i)} with $a_1 = a_2 = m + 1/2$, and the condition $a_1 + a_2\leq1$ is equivalent to $m\leq0$. The log-concavity assertion follows similarly from Theorem~\ref{thm:NCG.sum.shape}{\rm(iv)}.

\section*{Funding}
\addcontentsline{toc}{section}{Funding}

Robert E.~Gaunt is funded by EPSRC grant EP/Y008650/1. Fr\'ed\'eric Ouimet is supported by the Natural Sciences and Engineering Research Council of Canada (NSERC) through Discovery Grant RGPIN-2026-04471 and Discovery Launch Supplement DGECR-2026-00449.

\section*{References}

\setlength{\bibsep}{0pt plus 0ex}


\begin{thebibliography}{99}
\addcontentsline{toc}{section}{References}

\bibitem{a98} An, M. Y. Logconcavity versus logconvexity: a complete characterization. \emph{J. Econ. Theory} $\mathbf{80}$ (1998), 350--369.

\bibitem{bb05} Bagnoli, M. and Bergstrom, T. Log-concave probability and its applications. \emph{Econ. Theory} $\mathbf{26}$ (2005), 445--469.

\bibitem{BariczPonnusamyVuorinen2011} Baricz, \'{A}, Ponnusamy, S. and Vuorinen, M. Functional inequalities for modified Bessel functions. \emph{Expo. Math.} $\mathbf{29}$ (2011), 399--414.

\bibitem{cl05} Casta\~{n}o-Mart\'inez, A. and L\'opez-Bl\'azquez, F. Distribution of a Sum of Weighted Noncentral Chi-Square Variables. \emph{TEST} $\mathbf{14}$ (2005), 397--415.

\bibitem{craig} Craig, C. C. On the Frequency Function of $xy$. \emph{Ann. Math. Stat.} $\mathbf{7}$ (1936), 1--15.

\bibitem{cui} Cui, G., Yu, X., Iommelli, S. and Kong, L. Exact Distribution for the Product of Two Correlated Gaussian Random Variables. \emph{IEEE Signal Process. Lett.} $\mathbf{23}$ (2016), 1662--1666.

\bibitem{fw98} Feldmann, A. and Whitt, W. Fitting mixtures of exponentials to long-tail distributions to analyze network performance models. \emph{Perform. Eval.} $\mathbf{31}$ (1998), 245--279.

\bibitem{FischerGauntSarantsev2025VGReview} Fischer, A., Gaunt, R. E. and Sarantsev, A. The Variance-Gamma Distribution: A Review. \emph{Statist. Sci.} $\mathbf{40}$ (2025), 235--258.

\bibitem{gaunt prod} Gaunt, R. E. A note on the distribution of the product of zero mean correlated normal random variables. \emph{Stat. Neerl.} $\mathbf{73}$ (2019), 176--179.

\bibitem{g26} Gaunt, R. E. On the product of correlated normal random variables and the noncentral chi-square difference distribution. \emph{Statist. Probab. Lett.} $\mathbf{227}$ (2026), Art.\ 110554.

\bibitem{gm21} Gaunt, R. E. and Merkle, M. On bounds for the mode and median of the generalized hyperbolic and related distributions. \emph{J. Math. Anal. Appl.} $\mathbf{493}$ (2021), Art.\ 124508.

\bibitem{MR4924268} Gaunt, R. E., Nadarajah, S. and Pog\'any, T. K. Infinite divisibility of the product of two correlated normal random variables and exact distribution of the sample mean. \emph{J. Math. Anal. Appl.} $\mathbf{552}$ (2025), Art.\ 129800.

\bibitem{GauntSutcliffe2026NCGamma} Gaunt, R. E. and Sutcliffe, H. L. The non-central gamma sum and difference distributions: exact distribution and asymptotic expansions. arXiv:2605.15386, 2026.

\bibitem{gz23} Gaunt, R. E. and Ye, Z. Asymptotic approximations for the distribution of the product of correlated normal random variables. \emph{J. Math. Anal. Appl.} $\mathbf{543}$ (2025), Art.\ 128987.

\bibitem{pog1} G\'orska, K., Horzela, A., Jankov Ma\v{s}irevi\'c, D. and Pog\'any, T. K. Observations on the McKay $I_\nu$ Bessel distribution. \emph{J. Math. Anal. Appl.} $\mathbf{516}$ (2022), Art.\ 126481.

\bibitem{HolmAlouini2004} Holm, H. and Alouini, M.-S. Sum and difference of two squared correlated Nakagami variates in connection with the McKay distribution. \emph{IEEE Trans. Commun.} $\mathbf{52}$ (2004), 1367--1376.

\bibitem{j26} Jones, M. C. On mixture relationships between central and non-central chi-squared difference distributions. \emph{Stat. Probabil. Lett.} $\mathbf{229}$ (2026), Art.\ 110588.

\bibitem{k15} Klar, B. A note on gamma difference distributions. \emph{J. Stat. Comput. Simul.} $\mathbf{85}$ (2015), 3708--3715.

\bibitem{kkp01} Kotz, S., Kozubowski, T. J. and Podg\'{o}rski, K. \emph{The Laplace Distribution and Generalizations: A Revisit with Applications to Communications, Economics, Engineering, and Finance.} Birkh\"auser, Boston, 2001.

\bibitem{kb96} Kn\"{u}sel, L. and Bablok, B. Computation of the noncentral gamma distribution. \emph{SIAM J. Sci. Comput.} $\mathbf{17}$ (1996), 1224--1231.

\bibitem{lopez} L\'opez, J. L. Asymptotic expansions of integrals: The term-by-term integration method. \emph{J. Comput. Appl. Math.} $\mathbf{102}$ (1999), 181--194.

\bibitem{m93} Mathai, A. M. On noncentral generalized Laplacianness of quadratic forms in normal variables. \emph{J. Multivar. Anal.} $\mathbf{45}$ (1993), 239--246.

\bibitem{m32} McKay, A. T. A Bessel function distribution. \emph{Biometrika} $\mathbf{24}$ (1932), 39--44.

\bibitem{np16} Nadarajah, S. and Pog\'{a}ny, T. K. On the distribution of the product of correlated normal random variables. \emph{C. R. Acad. Sci. Paris, Ser. I} $\mathbf{354}$ (2016), 201--204.

\bibitem{MR2723248} Olver, F. W. J., Lozier, D. W., Boisvert, R. F. and Clark, C. W., eds. \emph{NIST Handbook of Mathematical Functions}. Cambridge University Press, New York, 2010.

\bibitem{p49} Patnaik, P. B. The non-central $\chi^2$- and $F$-distributions and their applications. \emph{Biometrika} $\mathbf{36}$ (1949), 202--232.

\bibitem{Prekopa1973} Pr\'{e}kopa, A. On logarithmic concave measures and functions. \emph{Acta Sci. Math. (Szeged)} $\mathbf{34}$ (1973), 335--343.

\bibitem{ss87} Shaked, M. and Shanthikumar, J. G. Characterization of some first passage times using log-concavity and log-convexity as aging notions. \emph{Probab. Eng. Inf. Sci.} $\mathbf{1}$ (1987), 279--291.

\bibitem{wb32} Wishart, J. and Bartlett, M. S. The distribution of second order moment statistics in a normal system. \emph{Proc. Cambridge Philos. Soc.} $\mathbf{28}$ (1932), 455--459.

\bibitem{y11} Yu, Y. The shape of the noncentral chi-square density. arXiv:1106.5241, 2011.

\bibitem{YuNVM2011} Yu, Y. On normal variance-mean mixtures. \emph{Statist. Probab. Lett.} $\mathbf{121}$ (2016), 45--50.

\end{thebibliography}
\end{document}